\documentclass[a4, 12pt]{amsart}
\usepackage[mathscr]{eucal}
\usepackage{amssymb}
\usepackage{latexsym}
\usepackage{amsthm}
\theoremstyle{plain}
\newtheorem{theorem}{Theorem}[section]

\newtheorem{remark}{Remark}[section]

\makeatletter
\@addtoreset{equation}{section}

\title[Chern conjecture]{Chern conjecture on minimal hypersurfaces}
\author{Qing-Ming Cheng, Guoxin Wei and Takuya Yamashiro}
\address{Qing-Ming Cheng \\ Department of Applied Mathematics, Faculty of Sciences,
Fukuoka  University, 814-0180, Fukuoka,  Japan, cheng@fukuoka-u.ac.jp}
\address{Guoxin Wei \\  School of Mathematical Sciences, South China Normal University,
510631, Guangzhou,  China, weiguoxin@tsinghua.org.cn}
\address{Takuya Yamashiro\\ Department of Applied Mathematics, Graduate School of Sciences,
Fukuoka  University, 814-0180, Fukuoka,  Japan, sd200002@cis.fukuoka-u.ac.jp}

\begin{document}
\maketitle

\begin{abstract}In this paper, we study $n$-dimensional complete
minimal hypersurfaces in a unit sphere.  We prove that an
$n$-dimensional complete  minimal hypersurface with constant scalar
curvature in a unit sphere with $f_3$ constant is isometric to  the  totally geodesic sphere or the
Clifford torus if $S\leq 1.8252 n-0.712898$, where $S$ denotes the
squared norm of the second fundamental form of this hypersurface.
\end{abstract}

\footnotetext{The first author was partially  supported by JSPS Grant-in-Aid for Scientific Research (B):  No.16H03937.
The second author was partly supported by grant No. 11771154 of NSFC, Guangdong Province Universities and Colleges Pearl River Scholar Funded Scheme (2018), Guangdong Natural Science Foundation Grant No.2019A1515011451.}

\maketitle
\section{ Introduction}
\noindent
As one  knows that it is   important  to investigate
compact minimal hypersurfaces in spheres. By computing  the Laplacian
of the squared norm $S$ of the second fundamental form of
minimal hypersurfaces in spheres, Simons in \cite{s}
proved that  for an $n$-dimensional compact
minimal hypersurface in a unit sphere $S^{n+1}(1)$, if $S\leq n$,
then $S\equiv 0$ or $S \equiv n$.  In the landmark papers  of
Chern, do Carmo and Kobayashi \cite{cdk}  and
Lawson \cite{l}, they proved that the Clifford torus $S^m(\sqrt{\frac mn})\times
S^{n-m}(\sqrt{\frac {n-m}n})$  for $1\leq m\leq n-1$  are the only compact minimal
hypersurfaces in $S^{n+1}(1)$ with  $S\equiv n$.
The following Chern conjecture is important and well-known:

\vskip2mm
\noindent
{\bf Chern  conjecture}.
For $n$-dimensional compact minimal  hypersurfaces in $S^{n+1}(1)$
with constant scalar curvature,  $S>n$, then  $S\geq 2n$.

\noindent
In 1982, Peng and Terng studied  the above Chern conjecture, they proved  that for $n$-dimensional compact minimal  hypersurfaces in $S^{n+1}(1)$
with constant scalar curvature,  $S>n$, then  $S\geq n+\dfrac1{12n}$. Furthermore, for $n=3$, they solved Chern conjecture affirmatively.
For $n\geq 4$, Yang and the first author (\cite{yc1}, \cite{yc2}, \cite{yc3}) made an important  breakthrough. They proved
if $S>n$, then $S\geq \dfrac{4n}3$. (cf. \cite{c}, \cite{dx}, \cite{gt}, \cite{gxxz}, \cite{pt}, \cite{pt1}, \cite{xx}).

\noindent
In \cite{cw},  Cheng and Wei have solved  Chern conjecture  for $n=4$
under the additional condition that  $f_3=\sum_i\lambda_i^3$ is constant, where $\lambda_i$
's are principal curvatures of $M^4$.
On the other hand, for $n\geq 4$,  Yang and the first author \cite{yc3} proved the following:
\begin{theorem}
Let $M^n$, $n\geq 4$,   be an  $n$-dimensional  compact minimal
hypersurface in $S^{n+1}(1)$ with constant scalar curvature. If
$f_3$ is constant,  $S=0$, or $S=n$ if $S\leq  n+\dfrac23n$.
\end{theorem}

\noindent
Our main purpose  in this paper is  to study  Chern  conjecture  under the condition that
$f_3$ is constant. We improve the result of Yang and Cheng \cite{yc3} under weaker topology.

\begin{theorem}\label{theorem 2}
Let $M^n$ $(n\geq 5)$ be an $n$-dimensional complete  minimal
hypersurface in $S^{n+1}(1)$ with constant scalar curvature. If
$f_3$ is constant and $S>n$, then
$$
S> 1.8252 n-0.712898.
$$
\end{theorem}
\begin{remark} In the above  theorem,  we only  assume that $M^n$
is complete.
\end{remark}

\section{Prelimenary}
\vskip2mm
\noindent
In  this paper, we assume that all manifolds are  smooth and
connected without boundary.
Let $M^n$ be an $n$-dimensional
hypersurface in $S^{n+1}(1)$. We choose a local orthonormal frame
$\{\vec{e}_1, \cdots, \vec{e}_{n}, \vec{e}_{n+1}\}$ and the dual
coframe $\{\omega_1, \cdots,$ $\omega_n$, $ \omega_{n+1}\}$ in
such a way that $\{\vec{e}_1, \cdots, \vec{e}_n\}$  is a local
orthonormal frame on $M^n$. Hence, we have
$$
\omega_{n+1}=0
$$
on $M^n$. Thus, one has
$$
\omega_{n+1,i}=\sum_{j=1}^nh_{ij}\omega_j, \ h_{ij}=h_{ji}.
$$
The mean curvature $H$ and the second fundamental form
$\vec{\alpha}$ of $M^n$ are defined, respectively, by
$$
H=\frac{1}{n}\sum_{i=1}^nh_{ii}, \
\vec{\alpha}=h_{ij}\omega_i\otimes\omega_j\vec{e}_{n+1}.
$$
If    $H$  is   zero in $M^n$, one
calls  that $M^n$ is   a  minimal hypersurface.
From the structure equations of $M^n$, Guass equations,
Codazzi equations and  Ricci formulas are given by
\begin{align*}
R_{ijkl}=(\delta_{ik}\delta_{jl}-\delta_{il}\delta_{jk})
+(h_{ik}h_{jl}-h_{il}h_{jk}),
\end{align*}
\begin{align*}
h_{ijk}=h_{ikj},
\end{align*}
\begin{align*}
h_{ijkl}-h_{ijlk}=\sum_mh_{im}R_{mjkl}+\sum_mh_{mj}R_{mikl}.
\end{align*}
where  $h_{ijk}=\nabla_kh_{ij}$ and $h_{ijkl}=\nabla_l\nabla_kh_{ij} $,
respectively.
For minimal hypersurfaces in $S^{n+1}(1)$ from
(2.1), we have
\begin{align*}
  r=n(n-1)-S,
\end{align*}
where $r$ and $S$ denote the scalar curvature and the squared norm
of the second fundamental form of $M^n$, respectively.
We define functions $f_3$ and $f_4$ by
\begin{align*}
f_3=\sum_{i,j,k=1}^nh_{ij}h_{jk}h_{ki} \ \ {\rm and} \ \
f_4=\sum_{i,j,k,l=1}^nh_{ij}h_{jk}h_{kl}h_{li}
\end{align*}
respectively.
 Then, we have, for minimal hypersurfaces,
  \begin{equation}
  \dfrac13\Delta f_3=(n-S)f_3+2C,
 \end{equation}
\begin{equation}
  \dfrac14\Delta f_4= (n-S)f_4+ (2A+B),
 \end{equation}
 where
 $$C=\sum_{i,j,k} \lambda_ih_{ijk}^2, \  \  A=\sum_{i,j,k}\lambda_i^2h_{ijk}^2, \ \
B=\sum_{i,j,k}\lambda_i\lambda_jh_{ijk}^2
$$
and $\lambda_i$'s are principal curvatures of $M^n$.
If the squared norm $S$ of the second fundamental
form is constant, we have
\begin{align*}
&\sum_ih_{ii}=\sum_i\lambda_i=0, \quad
S=\sum_{i,j}h_{ij}^2=\sum_i\lambda_i^2,\tag{2.7}\\
&\sum_{i,j,k}h_{ijk}^2= S(S-n),\ \ \sum_ph_{ijpp}=(n-S)h_{ij},\tag
{2.8}\\
&h_{ijij}-h_{jiji}=(\lambda_i-\lambda_j)(1+\lambda_i\lambda_j).\tag{2.9}
\end{align*}
By a direct computation, we have
\begin{align*}
\sum_{i,j,k,l}h_{ijkl}^2 +S(S-n)(2n+3-S) +3(2B-A)=0,\tag{2.10}
\end{align*}

\section{A proof of Theorem 1.2}

\vskip3mm
\noindent
Defining
\begin{align*}
  u_{ijkl}:=\frac 14(h_{ijkl}+h_{jkli}+h_{klij}+h_{lijk} ),
\end{align*}
we have
\begin{align*}
  \sum_{i,j,k,l}h_{ijkl}^2 = \sum_{i,j,k,l}u_{ijkl}^2 +\frac
  32\left[
  Sf_4-f_3^2-2S^2+nS\right],
\end{align*}
Putting $S-n=tS$ and defining
\begin{align*}
Sf = Sf_4-f_3^2-\frac {S^3}n,\tag{2.13}
\end{align*}
The following formulas can be found in \cite{yc2}
   \begin{equation}\label{eq:011} A-B \le \frac 13(\lambda_1-\lambda_2)^2tS^2,
\end{equation}
where $\lambda_1=\max_i\lambda_i$ and $\lambda_2=\min_i\lambda_i$.
Taking the orthonormal frame $\{\vec e_1, \vec e_2, \cdots, \vec e_n\}$ at each point such that
$$
h_{ij}=\lambda_i\delta_{ij},
$$
 we have
\begin{equation}\label{eq:12-23-5}
Sf+\frac{S^3}{n} \equiv Sf_4-\bigl(f_3\bigl)^2=S \sum_i\lambda_i^4-\left(
\sum_i\lambda_i^3\right)^2=\dfrac12\sum(\lambda_i-\lambda_j)^2\lambda_i^2\lambda_j^2.
\end{equation}
Since $S$ and $f_3$ are constant, we have
$$
\sum_i\lambda_ih_{iik}=0, \ \  \sum_i\lambda_i^2h_{iik}=0, \ \text{\rm for any  } \ k
$$
and
$$
\sum_i\lambda_ih_{iikl}= -\sum_{i,j}h_{ijk}h_{ijl}, \ \  \sum_i\lambda_i^2h_{iikl}=
-2\sum_{i,j}\lambda_ih_{ijk}h_{ijl}, \ \text{\rm for any  } \ k, l.
$$
Hence, we have
$$
\sum_{i,j,}h_{iijj} \lambda_i\lambda_j=-C, \ \ \sum_{i,j}h_{iijj}\lambda_i^2\lambda_j=-2B,  \ \
\sum_{i,j }h_{iijj}\lambda_i \lambda_j^2=-A
$$
Defining
$$
a_{ij}=\sum_mh_{im}h_{mj}-yh_{ij}-\dfrac Sn\delta_{ij},\ \ \ \ i, j=1,2,\cdots, n
$$
with $Sy=f_3$ we have
\begin{equation}
\begin{aligned}
&\sum_{i,j=1}^na_{ij}^2=f, \ \ \sum_{i,j=1}^na_{ij}h_{ij}=0, \ \ \sum_{i,j=1}^na_{ij}\delta_{ij}=0.
\end{aligned}
\end{equation}
For $ \forall \ \alpha, \beta, \gamma\in  \mathbb R$,

$$
\sum_{i,j,k,l}\biggl\{u_{ijkl}+\alpha (a_{ij}h_{kl}+h_{ij}a_{kl})+\beta h_{ij}h_{kl}+\gamma (h_{ij}\delta_{kl}+\delta_{ij} h_{kl})\biggl\}^2\geq 0,\
$$
Because of
\begin{equation}
\begin{aligned}
&\sum_{i,j,k,l}u_{ijkl}h_{ij}\delta_{kl}=-\dfrac{t}{2}S^2, \  \  \sum_{i,j,k,l}u_{ijkl}h_{ij}h_{kl}=-C,\\
&\sum_{i,j,k,l}u_{ijkl}a_{ij}h_{kl}=-B-\dfrac12A +yC+\dfrac{t}{2(1-t)}S^2,
\end{aligned}
\end{equation}
we obtain, from (3.3) and (3.4),
\begin{equation}
\begin{aligned}
\sum_{i,j,k,l}u_{ijkl}^2
\geq 2\alpha(2B+A-2yC-\dfrac{t}{1-t}S^2)-2\alpha^2Sf+\dfrac{C^2}{S^2}+\dfrac{t}{2(1-t)}tS^2
\end{aligned}
\end{equation}
by taking $\beta =\dfrac {C}{S^2}$ and $\gamma=\dfrac{t}{2(1-t)}$.
Since $f_3$ is constant, we have
$$
tSf_3=2C
$$
and
\begin{equation}
\begin{aligned}
&Sf+\dfrac{t}{1-t}S^2=Sf_4-f_3^2- S^2 =A-2B.
\end{aligned}
\end{equation}
From
\begin{equation}
\sum_{i,j,k,l}h_{ijkl}^2=S(S-n)(S-2n-3)+3(A-2B)
\end{equation}
and
$$
  \sum_{i,j,k,l}h_{ijkl}^2 = \sum_{i,j,k,l}u_{ijkl}^2 +\frac
  32\left[
  Sf_4-f_3^2-2S^2+nS\right],
  $$
  we have
\begin{equation}
S(S-n)(S-2n)=\sum_{i,j,k,l}u_{ijkl}^2+\frac 32S(S-n)-\dfrac32(A-2B).
\end{equation}
\begin{remark} If one  can prove
$$
\sum_{i,j,k,l}u_{ijkl}^2+\frac 32S(S-n)-\dfrac32(A-2B)\geq 0,
$$
Chern conjecture will be solved under condition that $f_3$ is constant.
\end{remark}

\noindent
Since $S$ is constant, we know that the Ricci curvature is bounded from below according to  the Gauss equations.
By applying  the Generalized Maximum Principle of  Omori \cite{o} and Yau \cite{y}  to function $f_4$,  we know that
there exists a sequence of $\{p_m\}_{m=1}^{\infty} \subset  M^n$  such that
$$
\lim_{m\to \infty}f(p_m)=\sup f_4, \lim_{m\to \infty}|\nabla f_4(p_m)|=0, \lim_{m\to \infty}\sup \Delta f_4(p_m)\leq 0
$$
Since $S$ is constant, we have that, for any $i, j, k, l$,
 $\{\lambda_i(p_m)\}$, $\{h_{ijk}(p_m)\}$ and $\{h_{ijkl}(p_m)\}$ are bounded sequences.
Thus, we can assume,  for any $i, j, k, l$,
 $$
 \lim_{m\to\infty}\lambda_i(p_m)=\bar{\lambda}_i, \ \  \lim_{m\to\infty} h_{ijk}(p_m) =\bar{h}_{ijk}
 \lim_{m\to\infty} h_{ijkl}(p_m)=\bar{h}_{ijkl}
 $$
 All of the following computations are made for $\bar{\lambda}_i$, $\bar{h}_{ijk}$ and $\bar{h}_{ijkl}$. For simple, we
 omit $\bar{ }$.
 From  (2.2), we have
$$
tSf_4\geq 2A+B.
$$
According to
$$
Sf_4-f_3^2-S^2=A-2B,
$$
we obtain
\begin{equation}\label{eq:12-23-1}
tf_3^2\geq  (2-t)A+(1+2t)B-tS^2.
\end{equation}
Since
\begin{equation}
\begin{aligned}
C^2&=(\sum_{i,j,k} \lambda_ih_{ijk}^2)^2=\dfrac19\bigl\{\sum_{i,j,k} (\lambda_i+\lambda_j+\lambda_k)h_{ijk}^2\bigl\}^2\\
&\leq  \dfrac13(A+2B)tS^2,
\end{aligned}
\end{equation}
Taking
\begin{equation}
\begin{aligned}
tz=\dfrac{17t^2-33t+24}{1-t}>0,
\end{aligned}
\end{equation}
we have
\begin{equation}\label{eq:12-23-2}
zC^2\leq \dfrac z3(A+2B)tS^2.
\end{equation}
From (3.5) and taking $\alpha =-\dfrac{3-4t}2$, we obtain
\begin{equation*}
\begin{aligned}
\ \ \ \ \ \ \ \ \ \ \ \ \ \ \ &\sum_{i,j,k,l}u_{ijkl}^2-\dfrac32(A-2B) +\frac 32S(S-n)\\
&\geq  2\alpha(2B+A)-(2\alpha^2+\dfrac32)(A-2B)+(-2\alpha t+(1+z)\dfrac{t^2}4)f^2_3-z\dfrac{C^2}{S^2}\\
&- 2\alpha\dfrac{t}{1-t}S^2+2\alpha^2\dfrac t{1-t}S^2 +\dfrac{t}{2(1-t)}tS^2+\dfrac 32tS^2\\
&\geq 2\alpha(2B+A)-(2\alpha^2+\dfrac32)(A-2B)\\
&+(-2\alpha +(1+z)\dfrac{t }4)\bigl\{(2-t)A+(1+2t)B\bigl\}-z\dfrac{t}3(A+2B)\\
&- 2\alpha\dfrac{t}{1-t}S^2+2\alpha^2\dfrac t{1-t}S^2 -(-2\alpha  +(1+z)\dfrac{t }4)tS^2+\dfrac{t}{2(1-t)}tS^2+\frac 32tS^2\\
&= \dfrac{4t^2-9t+3}{3(1-t)}(A-B)-2t\dfrac{t}{1-t}S^2,
\end{aligned}
\end{equation*}
that is,
\begin{equation}\label{eq:2-23-3}
\begin{aligned}
&S(S-n)(S-2n)\geq \dfrac{4t^2-9t+3}{3(1-t)}(A-B)-2t\dfrac{t}{1-t}S^2.
\end{aligned}
\end{equation}
If  $t\leq \dfrac{9-\sqrt {33}}8$, we have
$4t^2-9t+3\geq 0$. From \eqref{eq:2-23-3}, we get
$$
t\geq\frac{1}{2}-\frac{t}{(1-t)S},
$$
$$
S\geq 2n-\dfrac{2t}{(1-t)}\geq 2n- \dfrac{\sqrt {33}-3}{2},
$$
then
$$t\geq\frac{n}{2(n+1)}>\frac{5}{12}>\dfrac{9-\sqrt {33}}8.$$
It is a contradiction.
Hence, we have
$$
t> \dfrac{9-\sqrt {33}}8.
$$
From $ A-B\leq \dfrac23tS^3$, we obtain from \eqref{eq:2-23-3}
\begin{equation*}
\begin{aligned}
&(2t-1)S\geq \dfrac{4t^2-9t+3}{3(1-t)}\dfrac{2}{3}S-\dfrac{2t}{1-t},
\end{aligned}
\end{equation*}
that is,
\begin{equation*}
\begin{aligned}
26t^2-45t+15\leq \dfrac{18t}{S}.
\end{aligned}
\end{equation*}
Thus, we obtain
\begin{equation}
\begin{aligned}
t-\dfrac{45-\sqrt{465}}{52}\geq -\dfrac{36t}{45+\sqrt{465}-52t}\dfrac1{S}.
\end{aligned}
\end{equation}
We conclude
\begin{equation}
\begin{aligned}
t\geq  \dfrac{45-\sqrt{465}}{52}-\dfrac9{26}(\dfrac{3\sqrt{465}}{31}-1)\dfrac1S.
\end{aligned}
\end{equation}
In fact, if not,  we have
$$
t < \dfrac{45-\sqrt{465}}{52}-\dfrac9{26}(\dfrac{3\sqrt{465}}{31}-1)\dfrac1S<\dfrac{45-\sqrt{465}}{52}.
$$
Hence,  we infer
 $$
-\dfrac{36t}{45+\sqrt{465}-52t}\dfrac1{S}>-\dfrac9{26}(\dfrac{3\sqrt{465}}{31}-1)\dfrac1S.
 $$
 We  conclude
\begin{equation*}
\begin{aligned}
t> \dfrac{45-\sqrt{465}}{52}-\dfrac9{26}(\dfrac{3\sqrt{465}}{31}-1)\dfrac1S,
\end{aligned}
\end{equation*}
which is a contradiction. Thus, (3.15) must hold.\\
Hence, we get
 \begin{equation*}
\begin{aligned}
S-n=tS> \dfrac{45-\sqrt{465}}{52}S-\dfrac9{26}(\dfrac{3\sqrt{465}}{31}-1),
\end{aligned}
\end{equation*}
that is,
 \begin{equation}\label{eq:12-28-1}
 \begin{aligned}
S>\dfrac{\sqrt{465}-7}{8}n-\dfrac9{4}(1-\dfrac{\sqrt{465}}{31})\approx 1.82048 n-0.684881.
\end{aligned}
\end{equation}
Furthermore, we give a better estimate on  $t$. In order to do it,
 for any  $i , j$, we have
$$
 -\lambda_i\lambda_j\leq\dfrac14(\lambda_i-\lambda_j)^2.
$$
Hence, we get,  for $\lambda_i\lambda_j\leq 0$ and  $\lambda_i\lambda_k\leq 0$,
\begin{equation}\label{eq:12-23-4}
(|\lambda_i\lambda_j|+|\lambda_i\lambda_k|)^3\leq 4(|\lambda_i\lambda_j|^3+|\lambda_i\lambda_k|^3)
\leq (\lambda_i-\lambda_j)^2\lambda_i^2\lambda_j^2+(\lambda_i-\lambda_k)^2\lambda_i^2\lambda_k^2.
\end{equation}
For three different $i , j, k$, we know that at least one of
$\lambda_i\lambda_j$,  $\lambda_i\lambda_k$ and $\lambda_j\lambda_k$ is non-negative.
Without loss of generality, we assume
$\lambda_j\lambda_k\geq0$
and $\lambda_i\lambda_j\leq0$, $\lambda_i\lambda_k\leq0$,
then we get from \eqref{eq:12-23-4} and \eqref{eq:12-23-5}
\begin{equation}\label{eq:12-23-6}
-\lambda_i\lambda_j-\lambda_i\lambda_k-\lambda_j\lambda_k\leq |\lambda_i\lambda_j|
+|\lambda_i\lambda_k|
\leq (Sf+\frac{S^3}n)^{\frac13},
\end{equation}
and
\begin{equation}\label{eq:12-23-7}
-2\lambda_i\lambda_j\leq 2(|\lambda_i\lambda_j|^3)^{\frac13}
\leq 2(\dfrac14(\lambda_i-\lambda_j)^2\lambda_i^2\lambda_j^2)^{\frac13}
\leq 2\bigl\{\dfrac{Sf+\frac{S^3}n}4\bigl\}^{\frac13}.
\end{equation}
Since
\begin{equation}
\begin{aligned}
3(A-B)&=
\sum_{i,j,k}(\lambda_i^2+\lambda_j^2+\lambda_k^2-\lambda_i\lambda_j-\lambda_i\lambda_k-\lambda_j\lambda_k)h_{ijk}^2\\
&=\sum_{i\neq j\neq k\neq i}(\lambda_i^2+\lambda_j^2+\lambda_k^2-\lambda_i\lambda_j-\lambda_i\lambda_k-\lambda_j\lambda_k)h_{ijk}^2\\
&+3\sum_{i\neq j}(\lambda_i^2+\lambda_j^2-2\lambda_i\lambda_j)h_{iij}^2,
\end{aligned}
\end{equation}
we conclude from \eqref{eq:12-23-6} and \eqref{eq:12-23-7},
\begin{equation}
\begin{aligned}
3(A-B)&\leq (S+2\biggl\{\dfrac{Sf+\frac{S^3}n}4\biggl\}^{\frac13})\sum_{i,j,k}h_{ijk}^2
=(S+2\biggl\{\dfrac{Sf+\frac{S^3}n}4\biggl\}^{\frac13})tS^2.
\\
\end{aligned}
\end{equation}
According to
\begin{equation}\label{eq:12-23-8}
S-n=tS,  \  \ (1-t)Sf\leq \dfrac13(A-B), \  \ A-B\leq \dfrac23tS^3
\end{equation}
we infer from \eqref{eq:12-23-8}
\begin{equation}\label{eq:12-23-9}
\begin{aligned}
3(A-B)&\leq (S+2\biggl\{\dfrac{A-B}{12(1-t)}+\dfrac{1}{4n}\biggl\}^{\frac13})tS^2.
\end{aligned}
\end{equation}
\\
%


\noindent
Thus, we can  assume that  $3(A-B)\leq a_k t S^3$.  We have from \eqref{eq:12-23-9} that  $3(A-B)\leq a_{k+1} t S^3$, where
\begin{equation}\label{eq:4-22-10}
a_1=2,\ \ \ \ \ \ \    a_{k+1}=1+2\biggl(\frac{1}{36}\frac{t}{1-t}a_k+\frac{1}{4 n}\biggl)^{\frac{1}{3}}.
\end{equation}

\noindent We next assume that $t<0.452115$ and $n\geq 6$, then we get from \eqref{eq:4-22-10} that
\begin{equation}
 a_{k+1}\leq1+2\biggl(\frac{1}{36}\frac{0.452115}{1-0.452115}a_k+\frac{1}{24}\biggl)^{\frac{1}{3}}.
\end{equation}
By a direct calculation, we know
\begin{equation}
a_7\leq 1.878415.
\end{equation}
Hence, we obtain
\begin{equation}\label{eq:12-23-10}
3(A-B)<1.878415 t S^2.
\end{equation}
From \eqref{eq:2-23-3} and \eqref{eq:12-23-10}, we have
\begin{equation}\label{eq:12-23-11}
(2t-1)S> \dfrac{4t^2-9t+3}{3(1-t)}\times \frac{1.878415 }{3} S-\dfrac{2t}{1-t}.
\end{equation}
Then, we get
\begin{equation*}
2.83485(t-0.452115)(t-1.26876)<\dfrac{2t}{S}.
\end{equation*}
Because of
$t<0.452115$,  we know
 \begin{equation*}
\frac{2t}{2.83485(t-1.26876)}\frac{1}{S}>-0.390586\frac{1}{S}.
\end{equation*}
Thus, we infer
\begin{equation*}
t>0.452115-0.390586\frac{1}{S},
\end{equation*}
then
 \begin{equation}\label{eq:12-28-2}
 S>1.8252n-0.712898.
 \end{equation}
  If $t\geq0.452115$, then we get
  \begin{equation}
 S\geq1.8252n.
 \end{equation}
From \eqref{eq:12-28-1} and \eqref{eq:12-28-2}, we know that if $n=5$,
 $$
 1.82048 n-0.684881>1.8252n-0.712898.
$$
We complete the proof of theorem \ref{theorem 2}.

\bibliographystyle{amsplain}

\begin{thebibliography}{99}
\bibitem{c1}
Chang S. P., On minimal hypersurfaces with constant
scalar curvature in $S^4$, J. Diff. Geom., \textbf{37}(1993),
523-534.

\bibitem {c} Cheng Q.-M., The rigidity of Clifford torus
$S^1(\sqrt{\frac1n})\times S^{n-1}(\sqrt{\frac{n-1}{n}})$,
Comment. Math. Helv., \textbf{71}(1996), 60-69.

\bibitem{ci} Cheng Q.-M. \& Ishikawa, S., A characterization
of the Clifford torus, Proc. Amer. Math. Soc., \textbf{127}(1999),
819-828.

\bibitem {cw} Cheng Q.-M. and Wei G.,  Chern problems on minimal
hypersurfaces, preprint.

\bibitem {cdk} Chern S. S. do Carmo M.  \& Kobayashi S.,
Minimal submanifolds of a sphere with second fundamental form of
constant length, Functioal Analysis and Related Fields,
Springer-Verlag,  Berlin, 1970, pp. 59-75

\bibitem {dgw}
Deng Q. Gu H. \& Wei Q.,  Closed Willmore minimal hypersurfaces with constant scalar curvature in $S^5(1)$  are isoparametric,
 Adv. Math., \textbf{314}(2017), 278-305.

\bibitem {dx}
Ding Q. \& Xin Y. L., On Chern's problem for rigidity of minimal hypersurfaces in the spheres, Adv. Math.,  \textbf{227} (2011), 131-145.


\bibitem {gt}
Ge J. Q. \& Tang Z. Z., Chern conjecture and isoparametric hypersurfaces, Differential
geometry, Adv. Lect. Math., \textbf{22}, International Press, Somerville, MA, 2012, 49-60.


\bibitem {gxxz}
Gu J. Xu H. Xu Z. \& Zhao E., A survey on rigidity problems in geometry
and topology of submanifolds, Proceedings of the 6th International Congress of Chinese
Mathematicians, Adv. Lect. Math., \textbf{37}, Higher Education Press and International Press,
Beijing-Boston, 2016, 79-99.

\bibitem {l}  Lawson H. B. Jr., Local rigidity theorems for
minimal hypersurfaces,  Ann. of  Math., \textbf{89} (1969),
167-179.

\bibitem {o} Omori, H., Isometric immersions of Riemannian manifolds,
J. Math. Soc. Japan  \textbf{19} (1967), 205-214.



\bibitem {pt}  Peng C. K. \& Terng C. L., Minimal
hypersurfaces of spheres with constant scalar curvature, in ``Seminar on minimal submanifolds'', Princeton Univ. Press,
Princeton, 1983, pp. 179-198.

\bibitem {pt1}  Peng C. K. \& Terng C. L., The scalar curvature of
minimal hypersurfaces in sphere, Math. Ann., \textbf{266}(1983),
105-113.
\bibitem {s} Simons J., Minimal varieties
in Riemannian manifolds,  Ann. of Math., \textbf{ 88}(1968),
62-105.

%
%



\bibitem {xx}
Xu H. \& Xu Z., On Chern's conjecture for minimal hypersurfaces and rigidity of self-shrinkers, J. Funct. Anal., \textbf{273} (2017), 3406-3425.

\bibitem {yc1} Yang H. C. \& Cheng Q.-M., A note on the pinching constant
of minimal hypersurfaces with constant scalar curvature in the unit sphere, Chinese Science Bull., \textbf{36}(1991), 1-6.

\bibitem {yc2} Yang H. C. \& Cheng Q.-M., An estimate of the pinching constant of
minimal hypersurfaces with constant scalar curvature in the unit sphere, Manuscripta Math., \textbf{84}(1994), 89-100.

\bibitem {yc3} Yang H. C. \& Cheng Q.-M., Chern's conjecture on
minimal hypersurfaces, Math. Z., \textbf{227}(1998), 377-390.

\bibitem {y} Yau, S. T., Harmonic functions on complete Riemannian manifolds,
Comm. Pure and Appl. Math., \textbf{28} (1975), 201-228.


\end{thebibliography}

\end {document}